\documentclass[11pt]{article} 
\usepackage{latexsym} 
\newcommand{\eproof}{\mbox{\ }\hfill $\Box$ \par \vskip 10pt}

\newtheorem{Theorem}{Theorem}[section] 
\newtheorem{lemma}[Theorem]{Lemma} 
\newtheorem{prop}[Theorem]{Proposition}

\newtheorem{corol}[Theorem]{Corollary}

\begin{document}

\title{On the stabilization of the elasticity system by the boundary}

\author{{\sc Moez Khenissi and Georgi Vodev}}

\date{} 
 
\maketitle 

\begin{abstract}
\noindent
 We obtain free of resonances regions for the elasticity system in the exterior of a strictly convex body in
 ${\bf R}^3$ with dissipative boundary conditions under some natural assumptions on the behaviour of the geodesics
 on the boundary. To do so, we use the properties of the parametrix of the Neumann operator constructed in \cite{kn:SV1}. 
 As a consequence, we obtain time decay estimates for the local
 energy of the solutions of the corresponding mixed boundary value problems.
\end{abstract}

\setcounter{section}{0}
\section{Introduction and statement of results}

Let ${\cal O}\subset {\bf R}^3$ be a strictly convex compact set with 
smooth boundary $\Gamma=\partial{\cal O}$ and denote by $\Omega={\bf R}^3
\setminus{\cal O}$ the exterior domain. Denote by $\Delta_e$ the elasticity
operator, which is a $3\times 3$ matrix-valued differential operator defined
by
$$\Delta_eu=\mu_0\Delta u+(\lambda_0+\mu_0)\nabla(\nabla\cdot u),$$
$u=^t(u_1,u_2,u_3)$. Here $\lambda_0$, $\mu_0$ are the Lam\'e constants
supposed to satisfy
$$\mu_0>0,\quad 3\lambda_0+2\mu_0>0.\eqno{(1.1)}$$
The Neumann boundary conditions for $\Delta_e$ are of the form
$$(Bu)_i|_\Gamma:=\sum_{j=1}^3\sigma_{ij}(u)\nu_j|_\Gamma=0,\quad i=1,2,3,
\eqno{(1.2)}$$
where 
$$\sigma_{ij}(u)=\lambda_0\nabla\cdot u\delta_{ij}+\mu_0\left(\frac{
\partial u_i}{\partial x_j}+\frac{
\partial u_j}{\partial x_i}\right)$$
is the stress tensor, $\nu$ is the outer unit normal to $\Gamma$. The purpose
of the present paper is to study the time decay properties of the
elasticity system in $\Omega$ with dissipative boundary conditions. More
precisely, we are going to study the following mixed boundary value problem
$$\left\{
\begin{array}{l}
(\partial_t^2-\Delta_e)u=0\quad \mbox{in}\quad(0,+\infty)\times\Omega,\\
Bu-iAu=0\quad \mbox{on}\quad(0,+\infty)\times\Gamma,\\
u(0)=f_1,\,\partial_tu(0)=f_2,
\end{array}
\right.
\eqno{(1.3)}
$$
where $A$ is a classical zero order 
$3\times 3$ matrix-valued pseudo-differential 
operator on $\Gamma$, independent of $t$ and satisfying the
properties $A=A^*$, $A\ge 0$. Moreover, we suppose that there
exist a non-empty compact set $\Gamma_0\subset\Gamma$ and a constant
$C>0$ so that we have
$$\left\langle Af,f\right\rangle_{L^2(\Gamma)}\ge C\left\|
f\right\|^2_{L^2(\Gamma_0)}.\eqno{(1.4)}$$

The large time behaviour of the solutions to (1.3) with $A\equiv 0$ is
well understood. Kawashita \cite{kn:Ka} showed that there is no uniform
local energy, while Stefanov-Vodev \cite{kn:SV1}, \cite{kn:SV2} 
proved the existence of 
infinitely many resonances converging polynomially fast to the real axis. 
The reason for this is the existence of surface waves (called
Rayleigh waves), that is, a propagation of singularities of the
solutions along the geodesics on $\Gamma$ with a speed $c_R>0$ strictly 
less than the two other speeds in $\Omega$. Therefore, a strictly convex
obstacle is trapping for the Neumann problem of the elasticity wave equation.
Note that for the Dirichlet problem it is non-trapping, and in particular
we have an exponential decay of the local energy similarly to the classical
wave equation (see \cite{kn:Y}). Comming back to the equation (1.3) with non-trivial $A$,
note that we still have a propagation of singularities of the
solutions along the geodesics on $\Gamma$ with a speed $c_R>0$. Therefore,
in order to be able to get a better decay of the local energy we need to
suppose that all geodesics meet the part on $\Gamma$ where the 
dissipative term is non-trivial. More precisely, we suppose that there
exist a non-empty open domain $\Gamma'_0\subset\Gamma_0$ and a constant
$T>0$ so that
$$\mbox{for every geodesics}\,\,\gamma\,\, \mbox{with}\,\,
 \gamma(0)\in \Gamma,\,\, 
\mbox{there exists}\,\, 0\le t\le T \,\,\mbox{such that}
\,\, \gamma(t)\in\Gamma'_0.\eqno{(1.5)}$$
The outgoing resolvent, $R(\lambda)$, corresponding to the problem
(1.3) is defined via the equation
$$\left\{
\begin{array}{l}
(\Delta_e+\lambda^2)R(\lambda)f=f\quad \mbox{in}\quad\Omega,\\
\left(B-iA\right)R(\lambda)f=0\quad \mbox{on}\quad\Gamma,\\
R(\lambda)f - \lambda-\mbox{outgoing}.
\end{array}
\right.
\eqno{(1.6)}
$$
Recall that ``$\lambda$-outgoing'' means that there exist
$a\gg 1$ and a compactly supported function $g$ so that
$$R(\lambda)f|_{|x|\ge a}=R_0(\lambda)g|_{|x|\ge a},$$
where $R_0(\lambda)$ is the outgoing free resolvent, i.e. 
$$R_0(\lambda)=(\Delta_e+\lambda^2)^{-1}\in{\cal L}(L^2)\quad
\mbox{for}\quad {\rm Im}\,\lambda<0.$$
Let $\chi\in C_0^\infty({\bf R}^3)$, $\chi=1$ on ${\cal O}$. In the same way
as in the case $A\equiv 0$ we have that the cutoff resolvent
$$R_\chi(\lambda):=\chi R(\lambda)\chi$$
extends meromorphically to the whole complex plane ${\bf C}$ with poles
in ${\rm Im}\,\lambda>0$ called resonances. One of our goals in the present
paper is to study the distributions of the resonances near the real axis
under the above assumptions. Our first result is the following

\begin{Theorem} Under the assumptions (1.1), (1.4) and (1.5), 
$R_\chi(\lambda)$ extends analytically to
$\{|{\rm Im}\,\lambda|\le C_1|\lambda|^{-1},\,|{\rm Re}\,\lambda|\ge C_2>0\}$
and satisfies there the estimate
$$\left\|R_\chi(\lambda)\right\|_{{\cal L}(L^2)}\le C'.
\eqno{(1.7)}$$
Moreover, under the assumption (1.1) only, there exists a constant 
$C>0$ so that $R_\chi(\lambda)$ is analytic in the region 
$$\{C\le{\rm Im}\,\lambda\le M\log|\lambda|,\,|{\rm Re}\,\lambda|\ge C_M\gg 1
\}\eqno{(1.8)}$$ for every $M\gg 1$. Furthermore, there are infinitely many resonances in
$\{0<{\rm Im}\,\lambda< C\}$.
\end{Theorem}

In the case $A\equiv 0$, Stefanov-Vodev \cite{kn:SV1} showed that there is
a free of resonances region of the form
$$\{C_N|\lambda|^{-N}\le{\rm Im}\,\lambda\le M\log|\lambda|,
\,|{\rm Re}\,\lambda|\ge C_M\gg 1\}$$
for every $M,N\gg 1$, while in $\{0<{\rm Im}\,\lambda\le C_N|\lambda|^{-N}\}$
there are infinitely many resonances (called Rayleigh resonances) due to the 
Rayleigh surface waves. Later on Sj\"ostrand-Vodev \cite{kn:SjV} proved that
the counting function of these resonances is $$\tau_2c_R^{-2}{\rm Vol}(\Gamma)
r^2+O(r),\quad r\gg 1,\eqno{(1.9)}$$
where $c_R>0$ is the speed of the Rayleigh waves and 
$$\tau_2=(2\pi)^{-2}{\rm Vol}(\{x\in{\bf R}^3:|x|=1\}).$$ 
We expect that the counting function of the resonances in $\{0<{\rm Im}\,\lambda<C\}$ in the general case (i.e. for non-trivial $A$)
satisfies (1.9) as well with possibly an worse bound for the remainder term. 

On the other hand, extending a previous 
result by Burq \cite{kn:B}
to the elastic system, Bellassoued \cite{kn:Be} obtained a free of resonances
region of the form $\{0<{\rm Im}\,\lambda\le e^{-C|\lambda|}\}$, $C>0$,
so the Rayleigh resonances are concentrated in a region of the form
$\{e^{-C|\lambda|}\le {\rm Im}\,\lambda\le C_N|\lambda|^{-N}\}$.
Moreover, if the boundary $\Gamma$ is analytic, Vodev \cite{kn:V1} improved
this region to $\{e^{-C|\lambda|}\le {\rm Im}\,\lambda\le e^{-C'|\lambda|}\}$.
The presence of a non-trivial dissipative term $A$, however, changes the
distribution of the resonances considerably. 

As a consequence of (1.7) we get a decay rate of the local energy of
the solutions to (1.3).

\begin{corol} Under the assumptions (1.1), (1.4) and (1.5), for every
$a\gg 1$, $m\ge 0$, there exists a constant $C=C(a,m)>0$ so that we have
(for $t\gg 1$)
$$\|\nabla_xu(t,\cdot)\|_{L^2(\Omega_a)}+
\|\partial_t u(t,\cdot)\|_{L^2(\Omega_a)}\le C\left(t^{-1}\log t\right)^m
\left(\|\nabla f_1\|_{H^m(\Omega)}+\|f_2\|_{H^m(\Omega)}\right),
\eqno{(1.10)}$$
where $\Omega_a:=\Omega\cap\{|x|\le a\}$ and ${\rm supp}\,f_j\subset
\Omega_a$, $j=1,2$.
\end{corol}

The fact that (1.7) implies (1.10) was proved in \cite{kn:PV} in the case
of a unitary group. In our case this can be done following the
approach developed in \cite{kn:L} (and also in \cite{kn:B}). Note that in the case $A\equiv 0$, Bellassoued \cite{kn:Be}
proved (1.10) with $t^{-1}\log t$ replaced by $(\log t)^{-1}$.

It turns out that if the dissipation on the boundary is stronger, 
we have a uniform exponential decay of the local energy. Indeed, consider
the following mixed boundary value problem
$$\left\{
\begin{array}{l}
(\partial_t^2-\Delta_e)u=0\quad \mbox{in}\quad(0,+\infty)\times\Omega,\\
Bu+A\partial_tu=0\quad \mbox{on}\quad(0,+\infty)\times\Gamma,\\
u(0)=f_1,\,\partial_tu(0)=f_2,
\end{array}
\right.
\eqno{(1.11)}
$$
where $A$ is as above. The outgoing resolvent, $\widetilde R(\lambda)$, 
corresponding to the problem (1.11) is defined via the equation
$$\left\{
\begin{array}{l}
(\Delta_e+\lambda^2)\widetilde R(\lambda)f=f\quad \mbox{in}\quad\Omega,\\
\left(B-i\lambda A\right)\widetilde R(\lambda)f=0\quad \mbox{on}\quad\Gamma,\\
\widetilde R(\lambda)f - \lambda-\mbox{outgoing}.
\end{array}
\right.
\eqno{(1.12)}
$$
We have the following

\begin{Theorem} Assume (1.1) and (1.4) fulfilled with $\Gamma_0=\Gamma$. 
Then,  
$\widetilde R_\chi(\lambda)$ extends analytically to
$\{|{\rm Im}\,\lambda|\le C_1,\,|{\rm Re}\,\lambda|\ge C_2>0\}$
and satisfies there the estimate
$$\left\|\widetilde R_\chi(\lambda)\right\|_{{\cal L}(L^2)}
\le C'|\lambda|^{-1}.\eqno{(1.13)}$$
\end{Theorem}

As a consequence of (1.13) we get an exponential decay of the local 
energy of the solutions to (1.11).

\begin{corol} Under the assumptions of Theorem 1.3, for every
$a\gg 1$, there exist constants $C=C(a)>0$, $\alpha>0$, so that we have
(for $t\gg 1$)
$$\|\nabla_xu(t,\cdot)\|_{L^2(\Omega_a)}+
\|\partial_t u(t,\cdot)\|_{L^2(\Omega_a)}\le Ce^{-\alpha t}
\left(\|\nabla f_1\|_{L^2(\Omega)}+\|f_2\|_{L^2(\Omega)}\right),
\eqno{(1.14)}$$
provided ${\rm supp}\,f_j\subset\Omega_a$, $j=1,2$.
\end{corol}

The fact that (1.13) implies (1.14) is more or less well known in the
case of unitary groups (e.g. see \cite{kn:V2}). In the case of semi-groups
the proof goes in the same way (see \cite{kn:K}).

It is worth noticing that an interior dissipation of the elastic
wave equation with Neumann boundary conditions does not improve
the decay of the local energy. Indeed, consider the
following mixed boundary value problem
$$\left\{
\begin{array}{l}
(\partial_t^2-\Delta_e+A(x)\partial_t)u=0\quad \mbox{in}
\quad(0,+\infty)\times\Omega,\\
Bu=0\quad \mbox{on}\quad(0,+\infty)\times\Gamma,\\
u(0)=f_1,\,\partial_tu(0)=f_2,
\end{array}
\right.
\eqno{(1.15)}
$$
where $A\in C_0^\infty(\Omega)$ is a $3\times 3$ matrix-valued function  
satisfying the properties $A=A^*$, $A\ge 0$. Then, the quasi-modes 
constructed in \cite{kn:SV1}, \cite{kn:SV2}, which are due to the existence
of the Rayleigh waves and hence supported in an arbitrary small neighbourhood
of the boundary, are also quasi-modes for the problem with non-trivial $A$.
Therefore, in the same way as in these papers one can show that there
exists an infinite sequence $\{\lambda_j\}$ with $0<{\rm Im}\,\lambda_j\le
C_N|\lambda_j|^{-N}$, $\forall N\gg 1$, so that the following
problem has a non-trivial solution:
$$\left\{
\begin{array}{l}
(\Delta_e-i\lambda_j A(x)+\lambda_j^2)v_j=0\quad \mbox{in}\quad\Omega,\\
Bv_j=0\quad \mbox{on}\quad\Gamma,\\
v_j - \lambda_j-\mbox{outgoing}.
\end{array}
\right.
\eqno{(1.16)}
$$
Note finally that the situation is completely different for the usual scalar-valued wave equation with dissipative boundary
conditions like those above. Indeed, in this case if the obstacle is non-trapping, the corresponding cut-off resolvent
extends analytically through the real axis to a strip and as a consequence we have an exponential decay of the local energy
without extra assumptions (e.g. see \cite{kn:A}). In other words, the behaviour of the cut-off resolvent and the local energy is the same as in the case of the self-adjoint problem with Neumann boundary conditions.

\section{Proof of Theorem 1.1}

It sufices to prove (1.7) for real $\lambda\gg 1$, only. Let $v\in L^2_{comp}
(\Omega)$ and let $u$ be the solution to the equation
$$\left\{
\begin{array}{l}
(\Delta_e+\lambda^2)u=v\quad \mbox{in}\quad\Omega,\\
\left(B-iA\right)u=0\quad \mbox{on}\quad\Gamma,\\
u - \lambda-\mbox{outgoing}.
\end{array}
\right.
\eqno{(2.1)}
$$
Clearly, (1.7) is equivalent to the estimate
$$\|u\|_{L^2(\Omega_a)}\le C_a\|v\|_{L^2(\Omega)},\quad\lambda\ge\lambda_0,
\eqno{(2.2)}$$
for every $a\gg 1$ with constants
$C_a,\lambda_0>0$ indpendent of $\lambda$. 
To prove (2.2) we need a priori estimates
of the solutions to the equation
$$\left\{
\begin{array}{l}
(\Delta_e+\lambda^2)u=v\quad \mbox{in}\quad\Omega,\\
u|_{\Gamma}=f,\,\lambda^{-1}Bu|_\Gamma=g,\\
u - \lambda-\mbox{outgoing}.
\end{array}
\right.
\eqno{(2.3)}
$$
where $v\in L^2_{comp}(\Omega)$. We have the following

\begin{prop} There exist constants $C,\lambda_0>0$ so that for
$\lambda\ge\lambda_0$ we have
$$\|u\|_{H^1(\Omega_a)}+\|g\|_{L^2(\Gamma)}\le 
C\lambda^{-1}\|v\|_{L^2(\Omega)}+C\|f\|_{H^1(\Gamma)}.\eqno{(2.4)}$$
Hereafter the Sobolev spaces $H^1$ are equipped with the semi-classical norm
(with a small parameter $\lambda^{-1}$).
\end{prop}

 {\it Proof.} In the case of the Euclidean Laplacian $\Delta$ the 
a priori estimate (2.4) is proved in \cite{kn:CPV1} (see Theorem 3.1).
In our case the proof goes in the same way, but we will sketch it for
the sake of completeness. Observe first that the solution to (2.3)
is of the form
$$u=G(\lambda)v+K(\lambda)f,$$
where $G(\lambda)v$ solves the problem
$$\left\{
\begin{array}{l}
(\Delta_e+\lambda^2)G(\lambda)v=v\quad \mbox{in}\quad\Omega,\\
G(\lambda)v|_{\Gamma}=0,\\
G(\lambda)v - \lambda-\mbox{outgoing},
\end{array}
\right.
\eqno{(2.5)}
$$
while $K(\lambda)f$ solves the problem
$$\left\{
\begin{array}{l}
(\Delta_e+\lambda^2)K(\lambda)f=0\quad \mbox{in}\quad\Omega,\\
K(\lambda)f|_{\Gamma}=f,\\
K(\lambda)f - \lambda-\mbox{outgoing}.
\end{array}
\right.
\eqno{(2.6)}
$$
Since the strictly convex obstacles are non-trapping
for the Dirichlet problem of the elastic wave equation (see \cite{kn:Y}), 
we have the estimate
$$\|G(\lambda)v\|_{H^1(\Omega_a)}\le C_a\lambda^{-1}\|v\|_{L^2(\Omega)},
\quad\lambda\ge\lambda_0.\eqno{(2.7)}$$
Thus, to prove (2.4) we need the estimate
$$\|K(\lambda)f\|_{H^1(\Omega_a)}\le C_a\|f\|_{H^1(\Gamma)},
\quad\lambda\ge\lambda_0.\eqno{(2.8)}$$
This in turn follows from the fact that, since the obstacle is strictly
convex, one can construct a parametrix of $K(\lambda)$ near the boundary,
which satisfies (2.8). More precisely, there exist a neighbourhood
$\Omega'\subset\Omega$ of $\Gamma$ and operators
$${\cal K}(\lambda)=O(1):H^1(\Gamma)\to H^1(\Omega'),\quad
{\cal R}(\lambda)=O(\lambda^{-\infty}):H^1(\Gamma)\to H^1(\Omega'),
\eqno{(2.9)}$$
solving the equation
$$\left\{
\begin{array}{l}
(\Delta_e+\lambda^2){\cal K}(\lambda)f={\cal R}(\lambda)f
\quad \mbox{in}\quad\Omega',\\
{\cal K}(\lambda)f|_{\Gamma}=f.
\end{array}
\right.
\eqno{(2.10)}
$$
Note that such operators are constructed in \cite{kn:SV1} (Section 2). Let
$\psi\in C^\infty(\overline\Omega)$, supp$\,\psi\subset\overline\Omega'$,
$\psi=1$ near $\Gamma$. We have
$$\left\{
\begin{array}{l}
(\Delta_e+\lambda^2)\psi{\cal K}(\lambda)f=[\Delta_e,\psi]{\cal K}(\lambda)f
+\psi{\cal R}(\lambda)f
\quad \mbox{in}\quad\Omega,\\
\psi{\cal K}(\lambda)f|_{\Gamma}=f,
\end{array}
\right.
$$
which leads to
$$\left\{
\begin{array}{l}
(\Delta_e+\lambda^2)\left(K(\lambda)f-\psi{\cal K}(\lambda)f\right)
=-[\Delta_e,\psi]{\cal K}(\lambda)f
-\psi{\cal R}(\lambda)f
\quad \mbox{in}\quad\Omega,\\
\left(K(\lambda)f-\psi{\cal K}(\lambda)f\right)|_{\Gamma}=0.
\end{array}
\right.
$$
Hence
$$K(\lambda)f=\psi{\cal K}(\lambda)f-G(\lambda)\left([\Delta_e,\psi]{\cal K}
(\lambda)f+\psi{\cal R}(\lambda)f\right).\eqno{(2.11)}$$
Thus, (2.8) follows from (2.11), (2.7) and (2.9). To complete the proof
of (2.4) we need to show that
$$\|g\|_{L^2(\Gamma)}\le 
C\lambda^{-1}\|v\|_{L^2(\Omega)}+
\|u\|_{H^1(\Omega_a)}+C\|f\|_{H^1(\Gamma)}.\eqno{(2.12)}$$
To this end, we write the operator $\Delta_e$ in normal coordinates $y=(y_1,y')\in {\bf R}^+\times\Gamma$
in a neighbourhood of the boundary. We have
$$\partial_{x_j}=\nu_j(y')\partial_{y_1}+\beta_j(y)\cdot\nabla_{y'},\quad j=1,2,3,$$
where $\nu(y')=(\nu_1(y'),\nu_2(y'),\nu_3(y'))$ is the unit normal at $y'\in\Gamma$. Hence
$$\Delta_e={\cal A}(y')\partial_{y_1}^2+{\cal Q}(y,\partial_{y'})+{\cal Q}_1(y,\partial_{y}),$$
where ${\cal Q}$ and ${\cal Q}_1$ are second and first order differential operators, respectively, while
${\cal A}(y')$ is a symmetric matrix-valued function defined by
$$\left({\cal A}(y')u\right)_k=\mu_0u_k+(\lambda_0+\mu_0)\nu_k\sum_{j=1}^3\nu_ju_j,\quad k=1,2,3.$$
It is easy to check that 
$$\det {\cal A}(y')=\mu_0^2(\lambda_0+2\mu_0)>0.$$
Set
$$E(y_1)=\left\langle({\cal Q}+\lambda^2)(\psi u)(y_1,\cdot),(\psi u)(y_1,\cdot)\right\rangle_{L^2}
+\left\langle{\cal A}\partial_{y_1}(\psi u)(y_1,\cdot),\partial_{y_1}(\psi u)(y_1,\cdot)\right\rangle_{L^2},$$
$\psi$ being the function above. We have
$$\frac{dE(y_1)}{dy_1}=\left\langle[\partial_{y_1},{\cal Q}](\psi u)(y_1,\cdot),(\psi u)(y_1,\cdot)\right\rangle_{L^2}$$ $$
+2{\rm Re}\,\left\langle({\cal Q}+\lambda^2)(\psi u)(y_1,\cdot),\partial_{y_1}(\psi u)(y_1,\cdot)\right\rangle_{L^2}
+2{\rm Re}\,\left\langle{\cal A}\partial_{y_1}^2(\psi u)(y_1,\cdot),\partial_{y_1}(\psi u)(y_1,\cdot)\right\rangle_{L^2}$$
 $$=\left\langle[\partial_{y_1},{\cal Q}](\psi u)(y_1,\cdot),(\psi u)(y_1,\cdot)\right\rangle_{L^2}
-2{\rm Re}\,\left\langle{\cal Q}_1(\psi u)(y_1,\cdot),\partial_{y_1}(\psi u)(y_1,\cdot)\right\rangle_{L^2}$$ $$
+2{\rm Re}\,\left\langle(\Delta_e+\lambda^2)(\psi u)(y_1,\cdot),\partial_{y_1}(\psi u)(y_1,\cdot)\right\rangle_{L^2}.$$
Hence
$$E(0)=-\int_0^\infty \frac{dE(y_1)}{dy_1}dy_1\le O(\lambda^2)\|\psi u\|^2_{H^1(\Omega)}+O(1)\|(\Delta_e+\lambda^2)(\psi u)\|^2_{L^2(\Omega)}.$$
On the other hand
$$\left\|\partial_{y_1}(\psi u)(0,\cdot)\right\|^2_{L^2(\Gamma)}\le CE(0)+O(\lambda^2)\|(\psi u)(0,\cdot)\|^2_{H^1(\Gamma)},$$
with a constant $C>0$. Combining these estimates we get
$$\lambda^{-1}\left\|\partial_{y_1}u(0,\cdot)\right\|_{L^2(\Gamma)}\le O(1)\|u(0,\cdot)\|_{H^1(\Gamma)}
+O(1)\| u\|_{H^1(\Omega_a)}+O(\lambda^{-1})\|(\Delta_e+\lambda^2) u\|_{L^2(\Omega)},$$
which clearly implies (2.12).
\eproof

Set $w=G(\lambda)v$, where $v$ is as in (2.1). 
If $u$ is the solution to (2.1), then the function $u-w$ solves the equation
$$\left\{
\begin{array}{l}
(\Delta_e+\lambda^2)(u-w)=0\quad \mbox{in}\quad\Omega,\\
\left(B-iA\right)(u-w)=-Bw\quad \mbox{on}\quad\Gamma,\\
(u-w) - \lambda-\mbox{outgoing}.
\end{array}
\right.
\eqno{(2.13)}
$$
Set $f=u|_\Gamma=(u-w)|_\Gamma$, $g=-\lambda^{-1}Bw|_\Gamma$. By (2.4),
$$\|g\|_{L^2(\Gamma)}\le C\lambda^{-1}\|v\|_{L^2(\Omega)},
\quad\lambda\ge\lambda_0.\eqno{(2.14)}$$
Furthermore, we have
$$\lambda^{-1}B(u-w)|_\Gamma=N(\lambda)f,\eqno{(2.15)}$$
where $N(\lambda):H^1(\Gamma)\to L^2(\Gamma)$ is the outgoing Neumann
operator. Thus, we get that the function $f$ satisfies the equation
$$\left(N(\lambda)-i\lambda^{-1}A\right)f=g\eqno{(2.16)}$$
with $g$ satisfying (2.14). It is easy to see that (2.2) follows from 
combining (2.4), (2.14) and the following

\begin{Theorem} Under the assumptions (1.1), (1.4) and (1.5), there exist
constants $C,\lambda_0>0$ so that the solution to (2.16) satisfies the
estimate
$$\|f\|_{H^1(\Gamma)}\le C\lambda\|g\|_{L^2(\Gamma)},\quad\lambda\ge
\lambda_0.\eqno{(2.17)}$$
\end{Theorem}

{\it Proof.} Since the outgoing Neumann operator satisfies
$$-{\rm Im}\,\langle N(\lambda)f,f\rangle_{L^2(\Gamma)}\ge 0,\eqno{(2.18)}$$
we obtain 
$$\langle Af,f\rangle_{L^2(\Gamma)}\le -\lambda\,{\rm Im}\,\langle 
g,f\rangle_{L^2(\Gamma)}\le \beta^{-2}\lambda^2\|g\|^2_{L^2(\Gamma)}+
\beta^2\|f\|^2_{L^2(\Gamma)},\eqno{(2.19)}$$
for every $\beta>0$. By (1.4) and (2.19),
$$\|f\|_{L^2(\Gamma_0)}\le C\beta^{-1}\lambda\|g\|_{L^2(\Gamma)}+
\beta\|f\|_{L^2(\Gamma)}.\eqno{(2.20)}$$
Now, using (1.5) together with the properties of the outgoing Neumann operator,
we will prove the estimate
$$\|f\|_{H^1(\Gamma)}\le C\lambda\|g\|_{L^2(\Gamma)}+
C\|f\|_{L^2(\Gamma_0)}.\eqno{(2.21)}$$
Clearly, (2.17) follows from (2.21) and (2.20) provided $\beta$ is taken
small enough.

To prove (2.21) we will make use of the properties of the parametrix,
${\cal N}(\lambda)$, of $N(\lambda)$ constructed in Section 3 of 
\cite{kn:SV1}. First of all, we have
$$\left\|N(\lambda)f-{\cal N}(\lambda)f\right\|_{L^2(\Gamma)}\le
O(\lambda^{-\infty})\|f\|_{L^2(\Gamma)}.\eqno{(2.22)}$$
Moreover, ${\cal N}(\lambda)$ is a $\lambda-\Psi DO$ with a characteristic
variety $\Sigma=\{\zeta\in T^*\Gamma:\|\zeta\|=c_R^{-1}\}$ belonging to
the elliptic region of the corresponding boundary value problem. In the
region $\{\zeta\in T^*\Gamma:\|\zeta\|>c_R^{-1}\}$ the operator
${\cal N}(\lambda)$ is an elliptic $\lambda-\Psi DO$ of class
$L^{1,0}_{0,0}(\Gamma)$ (hereafter we use the same notations as in the
appendix of \cite{kn:SV1}), while in the
region $\{\zeta\in T^*\Gamma:\|\zeta\|<c_R^{-1}\}$ it is hypoelliptic.
Clearly, so is the operator ${\cal N}(\lambda)-i\lambda^{-1}A$.
Therefore, if $\chi\in C_0^\infty(T^*\Gamma)$, $\chi=1$ on 
$\{\zeta\in T^*\Gamma:\left|\|\zeta\|-c_R^{-1}\right|\le\epsilon\}$,
$\chi=0$ on 
$\{\zeta\in T^*\Gamma:\left|\|\zeta\|-c_R^{-1}\right|\ge 2\epsilon\}$,
$0<\epsilon\ll 1$, we have
$$\|{\rm Op}_\lambda(1-\chi)f\|_{H^1(\Gamma)}\le O(\lambda^{2/3})
\left\|\left({\cal N}(\lambda)-i\lambda^{-1}A\right)f\right\|_{L^2(\Gamma)}
+O(\lambda^{-\infty})\|f\|_{H^1(\Gamma)}$$ $$
\le O(\lambda^{2/3})\|g\|_{L^2(\Gamma)}
+O(\lambda^{-\infty})\|f\|_{H^1(\Gamma)}.\eqno{(2.23)}$$
On the other hand, near $\Sigma$ the operator ${\cal N}(\lambda)$ 
is a $\lambda-\Psi DO$  of class $L^{0,0}_{0,0}(\Gamma)$, whose 
principal symbol is a symmetric $3\times 3$ matrix-valued function 
with eigenvalues $a_1(\zeta)=\widetilde a_1(\zeta)(c_R\|\zeta\|-1)$, 
$\widetilde a_1(\zeta)> 0$, $a_2(\zeta)>0$, $a_3(\zeta)> 0$
near $\Sigma$. It is shown in \cite{kn:S} (Theorem 3.1) that 
there exist elliptic $\lambda-\Psi DOs$, 
$U(\lambda)$ and $V(\lambda)$,  
of class $L^{0,0}_{0,0}(\Gamma)$, so that we have
$$U(\lambda)^*{\cal N}(\lambda){\rm Op}_\lambda(\widetilde\chi)U(\lambda)=
\left(
\begin{array}{lr}
-\lambda^{-2}c_R^{2}\Delta_\Gamma-\lambda^{-1}a_0-1& 0\\
0&V(\lambda)
\end{array}
\right){\rm Op}_\lambda(\widetilde\chi_1)
 \quad{\rm mod}\quad L^{0,-2}_{0,0}(\Gamma),
\eqno{(2.24)}$$
where $-\Delta_\Gamma$ denotes the positive Laplace-Beltrami operator
on $\Gamma$, $a_0$ is a classical (independent of $\lambda$) zero order
$\Psi DO$ on $\Gamma$ with a real-valued principal symbol, 
and $\widetilde\chi,\widetilde\chi_1\in C_0^\infty(T^*\Gamma)$, $\widetilde
\chi_1=1$ on 
$\{\zeta\in T^*\Gamma:\left|\|\zeta\|-c_R^{-1}\right|\le 3\epsilon\}$,
$\widetilde\chi_1=0$ on 
$\{\zeta\in T^*\Gamma:\left|\|\zeta\|-c_R^{-1}\right|\ge 4\epsilon\}$, $\widetilde
\chi=1$ on supp$\,\widetilde\chi_1$.
In fact in \cite{kn:S} a better diagonalization of ${\cal N}(\lambda)$
 near $\Sigma$ is carried out, but for our purposes (2.24) will suffice. 
Now the function $\widetilde f=
U(\lambda)^{-1}{\rm Op}_\lambda(\chi)f$ satisfies 
$$\left(
\begin{array}{lr}
-\lambda^{-2}c_R^{2}\Delta_\Gamma-\lambda^{-1}a_0-1& 0\\
0&V(\lambda)
\end{array}
\right)\widetilde f-i\lambda^{-1}\widetilde A\widetilde f=
\widetilde g,\eqno{(2.25)}$$
where $\widetilde A=U^*AU$ is a $\lambda-\Psi DO$ of class
$L^{0,0}_{0,0}(\Gamma)$ with a principal symbol satisfying $\sigma_p
(\widetilde A)\ge 0$, $\sigma_p
(\widetilde A)^*=\sigma_p
(\widetilde A)$, and $\widetilde g$ satisfies
$$\|\widetilde g\|_{L^2(\Gamma)}\le \|g\|_{L^2(\Gamma)}
+O(\lambda^{-\infty})\|f\|_{H^1(\Gamma)}.\eqno{(2.26)}$$
Writing $\widetilde f=(\widetilde f_1,\widetilde f_2,\widetilde f_3)$,
$\widetilde g=(\widetilde g_1,\widetilde g_2,\widetilde g_3)$,
we reduce (2.25) to 
$$\left(-\lambda^{-2}c_R^{2}\Delta_\Gamma-1-i\lambda^{-1}b_1\right)
\widetilde f_1=\widetilde g_1+\lambda^{-1}\left(b_2\widetilde f_2+
b_3\widetilde f_3\right),\eqno{(2.27)}$$
$$\widetilde V(\lambda)
\left(
\begin{array}{l}
\widetilde f_2\\
\widetilde f_3
\end{array}
\right)
=
\left(
\begin{array}{l}
\widetilde g_2\\
\widetilde g_3
\end{array}
\right)
+\lambda^{-1}
\left(
\begin{array}{l}
c_2\widetilde f_1\\
c_3\widetilde f_1
\end{array}
\right),\eqno{(2.28)}$$
where $b_j$, $c_j$ are scalar-valued $\lambda-\Psi DOs$ of class
$L^{0,0}_{0,0}(\Gamma)$, the principal symbol of $b_1$ satisfying
${\rm Re}\,\sigma_p(b_1)\ge 0$, while $\widetilde V(\lambda)$ is a 
$2\times 2$ matrix-valued elliptic $\lambda-\Psi DO$ of class
$L^{0,0}_{0,0}(\Gamma)$. Thus, the inverse $\widetilde V(\lambda)^{-1}$ 
is again a $2\times 2$ matrix-valued elliptic $\lambda-\Psi DO$ of class
$L^{0,0}_{0,0}(\Gamma)$, so we can solve the equation (2.28).
In particular, we obtain
$$\|\widetilde f_2\|_{L^2(\Gamma)}+\|\widetilde f_3\|_{L^2(\Gamma)}
\le \|\widetilde g_2\|_{L^2(\Gamma)}+\|\widetilde g_3\|_{L^2(\Gamma)}
+O(\lambda^{-1})\|\widetilde f_1\|_{L^2(\Gamma)}.\eqno{(2.29)}$$
Furthermore, we conclude that the function $\widetilde f_1$
solves an equation of the form
$$\left(-\lambda^{-2}c_R^{2}\Delta_\Gamma-1-i\lambda^{-1}b\right)
\widetilde f_1=h,\eqno{(2.30)}$$
 where $b$ is a scalar-valued $\lambda-\Psi DO$ of class
$L^{0,0}_{0,0}(\Gamma)$ with a principal symbol satisfying
${\rm Re}\,\sigma_p(b)\ge 0$, and $h$ satisfies
$$\|h\|_{L^2(\Gamma)}\le \|g\|_{L^2(\Gamma)}
+O(\lambda^{-\infty})\|f\|_{H^1(\Gamma)}.\eqno{(2.31)}$$
We are going to show that the assumption (1.5) leads to the estimate
$$\|\widetilde f_1\|_{L^2(\Gamma)}\le C\lambda\|h\|_{L^2(\Gamma)}+
C\|\widetilde f_1\|_{L^2(\Gamma'_0)}.\eqno{(2.32)}$$
Before doing so, observe that (2.32) implies (2.21). Indeed, since
$$\|\widetilde f_1\|_{L^2(\Gamma'_0)}\le C\|f\|_{L^2(\Gamma_0)}
+O(\lambda^{-\infty})\|f\|_{L^2(\Gamma)},$$
we deduce from (2.29), (2.31) and (2.32) that
$$\|{\rm Op}_\lambda(\chi)f\|_{L^2(\Gamma)}\le C\lambda\|g\|_{L^2(\Gamma)}
+C\|f\|_{L^2(\Gamma_0)}
+O(\lambda^{-\infty})\|f\|_{H^1(\Gamma)}.\eqno{(2.33)}$$
Thus, (2.21) follows from (2.23) and (2.33).

The fact that (1.5) implies (2.32) can be derived from the more general results
of \cite{kn:BLR}, but we will give here a simpler proof following \cite{kn:V3}
where this is carried out in the case $b\equiv 0$ (see Theorem 2.3 of 
\cite{kn:V3}). Denote by $r_0(x,\xi)$, $(x,\xi)\in T^*\Gamma$, the 
principal symbol of the operator $-c_R^{2}\Delta_\Gamma$, so that we have
$\Sigma=\{(x,\xi)\in T^*\Gamma:r_0(x,\xi)=1\}$. Recall that the 
bicharacteristic flow $\Phi(t):T^*\Gamma\to T^*\Gamma$, $t\in{\bf R}$, 
associated
to the Hamiltonian $r_0(x,\xi)$ is defined by $\Phi(t)(x^0,\xi^0):=
(x(t),\xi(t))$, where the pair $(x(t),\xi(t))$ solves the Hamilton 
equation
$$\frac{\partial x(t)}{\partial t}=\frac{\partial r_0(x,\xi)}{\partial\xi},\,
\,\frac{\partial \xi(t)}{\partial t}=-\frac{\partial r_0(x,\xi)}{\partial x},\,
\,\, x(0)=x^0,\,\xi(0)=\xi^0.$$
Fix a point $\zeta^0=(x^0,\xi^0)\in\Sigma$ and choose a real-valued function 
$p(x,\xi)\in C_0^\infty(T^*\Gamma)$, $0\le p\le 1$, such that $p=1$ in
a neighbourhood of $\zeta^0$ and $p=0$ outside a biger neighbourhood. Given
a $t\in{\bf R}$, define the function $p_t(x,\xi)\in C_0^\infty(T^*\Gamma)$
by $p_t(x,\xi)=p(\Phi(-t)(x,\xi))$. By a microlocal partition of the unity
in a neighbourhood of $\Sigma$, it is easy to see that (2.32) follows from
(1.5) and the following

\begin{lemma} For every $T>0$ there exist positive constants $C=C(T)$
and $\lambda_0=\lambda_0(T)$ so that the solutions to (2.30) satisfy the
estimate
$$\left\|p(x,{\cal D}_x)\widetilde f_1\right\|_{L^2(\Gamma)}\le
\left\|p_t(x,{\cal D}_x)\widetilde f_1\right\|_{L^2(\Gamma)}+2T\lambda
\left\|h\right\|_{L^2(\Gamma)}+C\lambda^{-1}\left\|\widetilde f_1
\right\|_{L^2(\Gamma)},\eqno{(2.34)}$$
for $0\le t\le T$, $\lambda\ge\lambda_0$. Hereafter we denote
${\cal D}_x:=\lambda^{-1}D_x$.
\end{lemma}

{\it Proof.} Set $P=-\lambda^{-2}c_R^{2}\Delta_\Gamma-1-i\lambda^{-1}b$.
Since
$$\partial_tp_t+\{r_0,p_t\}=0,$$
the operator
$$Q_t:=\lambda \partial_tp_t(x,{\cal D}_x)+i\lambda^2[P,p_t(x,{\cal D}_x)]$$
is a zero order $\lambda-\Psi DO$, and hence uniformly bounded on 
$L^2(\Gamma)$. Moreover, the fact that the principal symbol of the 
operator $b$ satisfies ${\rm Re}\,\sigma_p(b)\ge 0$ implies
$$-{\rm Re}\,\langle bf,f\rangle_{L^2(\Gamma)}\le O(\lambda^{-1})\|f
\|_{L^2(\Gamma)}^2,\quad\forall f\in L^2(\Gamma).$$
Therefore, using the identity
$$\frac{1}{2}\frac{d}{dt}\left\|p_t(x,{\cal D}_x)\widetilde f_1
\right\|_{L^2(\Gamma)}^2={\rm Re}\,\left\langle 
\partial_tp_t(x,{\cal D}_x)\widetilde f_1,
p_t(x,{\cal D}_x)\widetilde f_1\right\rangle_{L^2(\Gamma)}$$
 $$=\lambda{\rm Re}\,\left\langle 
[P,p_t(x,{\cal D}_x)]\widetilde f_1,
p_t(x,{\cal D}_x)\widetilde f_1\right\rangle_{L^2(\Gamma)}+
\lambda^{-1}{\rm Re}\,\left\langle 
Q_t\widetilde f_1,
p_t(x,{\cal D}_x)\widetilde f_1\right\rangle_{L^2(\Gamma)}$$ 
$$=-\lambda{\rm Re}\,\left\langle 
p_t(x,{\cal D}_x)P\widetilde f_1,
p_t(x,{\cal D}_x)\widetilde f_1\right\rangle_{L^2(\Gamma)}+
\lambda^{-1}{\rm Re}\,\left\langle 
Q_t\widetilde f_1,
p_t(x,{\cal D}_x)\widetilde f_1\right\rangle_{L^2(\Gamma)}$$ $$
-{\rm Re}\,\left\langle 
 bp_t(x,{\cal D}_x)\widetilde f_1,
p_t(x,{\cal D}_x)\widetilde f_1\right\rangle_{L^2(\Gamma)},$$
we obtain
$$\left|\frac{d}{dt}\left\|p_t(x,{\cal D}_x)\widetilde f_1
\right\|_{L^2(\Gamma)}\right|\le 2\lambda
\left\|P\widetilde f_1\right\|_{L^2(\Gamma)}+O(\lambda^{-1})
\left\|\widetilde f_1\right\|_{L^2(\Gamma)}.\eqno{(2.35)}$$
By (2.35),
$$\left\|p(x,{\cal D}_x)\widetilde f_1
\right\|_{L^2(\Gamma)}=\left\|p_t(x,{\cal D}_x)\widetilde f_1
\right\|_{L^2(\Gamma)}-\int_0^t\frac{d}{d\tau}\left\|p_\tau(x,{\cal D}_x)
\widetilde f_1\right\|_{L^2(\Gamma)}d\tau$$
 $$\le \left\|p_t(x,{\cal D}_x)\widetilde f_1\right\|_{L^2(\Gamma)}
+2t\lambda
\left\|P\widetilde f_1\right\|_{L^2(\Gamma)}+O(\lambda^{-1})
\left\|\widetilde f_1\right\|_{L^2(\Gamma)}.$$
\eproof

Clearly, the free of resonances region follows from the following

\begin{prop} Under the assumption (1.1), for $\lambda$ belonging to the
region (1.8) with a suitably chosen constant $C>0$,  
the solution to (2.16) satisfies the estimate
$$\|f\|_{L^2(\Gamma)}\le \frac{C'|\lambda|}{{\rm Im}\,\lambda}
\|g\|_{L^2(\Gamma)}.\eqno{(2.36)}$$
\end{prop}

{\it Proof.} Without loss of generality we may suppose $\lambda_1:=
{\rm Re}\,\lambda\gg 1$. It is shown in \cite{kn:SV1} that, for
$\lambda$ belonging to the region $\Lambda_M=\{0\le {\rm Im}\,\lambda\le
M\log|\lambda|,\,{\rm Re}\,\lambda\ge C_M\gg 1\}$, the Neumann operator
has a parametrix ${\cal N}(\lambda)$ which is a $\lambda_1-\Psi DO$ 
with a characteristic variety $\Sigma$, depending on a parameter 
$\lambda_1^{-1}{\rm Im}\,\lambda\ll 1$. In particular, (2.22) still holds
with $O(\lambda^{-\infty})$ replaced by $O(\lambda_1^{-\infty})$. In what
follows we will keep the same notations as in the proof of Theorem 2.2 above.
In fact, much of the analysis still works with ${\rm Op}_\lambda$ replaced
by ${\rm Op}_{\lambda_1}$. For example, we have the following analogue of
(2.23)
$$\|{\rm Op}_{\lambda_1}(1-\chi)f\|_{H^1(\Gamma)}
\le O(|\lambda|^{2/3})\|g\|_{L^2(\Gamma)}
+O(\lambda_1^{-\infty})\|f\|_{H^1(\Gamma)},\quad \lambda\in\Lambda_M.
\eqno{(2.37)}$$
We still have (2.29) with $O(\lambda^{-1})$ replaced by $O(\lambda_1^{-1})$
as well as (2.30) with $h$ satisfying (2.31) with $O(\lambda^{-\infty})$ 
replaced by $O(\lambda_1^{-\infty})$. Thus, we get
$${\rm Im}\,\lambda^2\|\widetilde f_1\|^2_{L^2(\Gamma)}=-{\rm Im}\,\left
\langle \lambda^2h,\widetilde f_1\right\rangle_{L^2(\Gamma)}+
{\rm Re}\,\left
\langle \lambda b\widetilde f_1,\widetilde f_1\right\rangle_{L^2(\Gamma)}$$
 $$\le |\lambda|^2\|h\|_{L^2(\Gamma)}\|\widetilde f_1\|_{L^2(\Gamma)}+
C|\lambda|\|\widetilde f_1\|^2_{L^2(\Gamma)},$$
and hence
$${\rm Im}\,\lambda\|\widetilde f_1\|_{L^2(\Gamma)}\le \left(2
{\rm Im}\,\lambda-C\right)\|\widetilde f_1\|_{L^2(\Gamma)}\le O(|\lambda|)
\|h\|_{L^2(\Gamma)},\eqno{(2.38)}$$
provided ${\rm Im}\,\lambda\ge C$, $\lambda\in\Lambda_M$. Combining
(2.38) with (2.29) and (2.31), we obtain
$$\|{\rm Op}_{\lambda_1}(\chi)f\|_{L^2(\Gamma)}
\le \frac{C'|\lambda|}{{\rm Im}\,\lambda}\|g\|_{L^2(\Gamma)}
+O(\lambda_1^{-\infty})\|f\|_{H^1(\Gamma)},\eqno{(2.39)}$$
for $\lambda$ belonging to the region (1.8). Now (2.36) follows from
(2.37) and (2.39).

To prove the existence of infinitely many resonances (i.e. poles of $(\lambda N(\lambda)-iA)^{-1}$)
 in $\{0<{\rm Im}\,\lambda< C\}$ we will proceed as in \cite{kn:SV1}. Without loss of generality we may suppose
 ${\rm Re}\,\lambda>0$. By (2.36), we have
 $$\left\|(\lambda N(\lambda)-iA)^{-1}\right\|_{L^2(\Gamma)\to L^2(\Gamma)}\le C(\log|\lambda|)^{-1},\quad\lambda\in l^\pm,\eqno{(2.40)}$$
 where $l^\pm:=\{\lambda\in {\bf C}:\pm{\rm Im}\,\lambda=\log{\rm Re}
\,\lambda,\,{\rm Re}\,\lambda\ge C'\}$ with some constant
 $C'\gg 1$. If we suppose that $(\lambda N(\lambda)-iA)^{-1}$ is analytic in $\{\lambda\in {\bf C}:0<{\rm Im}\,\lambda<C,\,{\rm Re}\,\lambda\ge C'\}$, so it is in $\{\lambda\in {\bf C}:|{\rm Im}\,\lambda|\le\log{\rm Re}\,\lambda,\,{\rm Re}\,\lambda\ge C'\}$. Then, by (2.40) together with the Fragm\`en-Lindel\"of principle we get
 $$\left\|(\lambda N(\lambda)-iA)^{-1}\right\|_{L^2(\Gamma)\to L^2(\Gamma)}\le C(\log\lambda)^{-1},\quad\lambda\in {\bf R},\,\lambda\ge C'.\eqno{(2.41)}$$
 On the other hand, it is shown in \cite{kn:SV1} that there exist quasi-modes $(f_j,k_j)\in L^2(\Gamma)\times{\bf R}$ such that
 $\|f_j\|_{L^2}=1$, $k_j\to +\infty$ and
 $$\|k_jN(k_j)f_j\|_{L^2}\le Const.$$
 Hence,
 $$\|(k_jN(k_j)-iA)f_j\|_{L^2}\le Const,$$
 which combined with (2.41) lead to 
 $$1=\|f_j\|_{L^2}\le Const (\log k_j)^{-1},$$
 which is impossible if we take $k_j$ large enough. Therefore, the operator-valued function $(\lambda N(\lambda)-iA)^{-1}$ cannot be  analytic in $\{\lambda\in {\bf C}:0<{\rm Im}\,\lambda<C,\,{\rm Re}\,\lambda\ge C'\}$.
\eproof

\section{ Proof of Theorem 1.3} 

Again, it suffices to prove (1.13) for real
$\lambda\gg 1$, only.  Let $v\in L^2_{comp}
(\Omega)$ and let $u$ be the solution to the equation
$$\left\{
\begin{array}{l}
(\Delta_e+\lambda^2)u=v\quad \mbox{in}\quad\Omega,\\
\left(B-i\lambda A\right)u=0\quad \mbox{on}\quad\Gamma,\\
u - \lambda-\mbox{outgoing}.
\end{array}
\right.
\eqno{(3.1)}
$$
Clearly, (1.13) is equivalent to the estimate
$$\|u\|_{L^2(\Omega_a)}\le C_a\lambda^{-1}
\|v\|_{L^2(\Omega)},\quad\lambda\ge\lambda_0.\eqno{(3.2)}$$
The function $f=u|_\Gamma$ solves the equation
$$\left(N(\lambda)-iA\right)f=g\eqno{(3.3)}$$
with $g$ satisfying
 $$\|g\|_{L^2(\Gamma)}\le C\lambda^{-1}\|v\|_{L^2(\Omega)},
\quad\lambda\ge\lambda_0.\eqno{(3.4)}$$
Thus, in view of (2.4), to prove (3.2) it suffices to show that
$$\|f\|_{H^1(\Gamma)}\le C\|g\|_{L^2(\Gamma)},\quad\lambda\ge
\lambda_0.\eqno{(3.5)}$$
Using (2.18) and (1.4) with $\Gamma_0=\Gamma$, we get
$$C\|f\|^2_{L^2(\Gamma)}\le\langle Af,f\rangle_{L^2(\Gamma)}
\le -{\rm Im}\,\langle 
g,f\rangle_{L^2(\Gamma)}\le \beta^{-2}\|g\|^2_{L^2(\Gamma)}+
\beta^2\|f\|^2_{L^2(\Gamma)},\eqno{(3.6)}$$
for every $\beta>0$. Taking $\beta$ small enough, we deduce from (3.6),
$$\|f\|_{L^2(\Gamma)}\le C\|g\|_{L^2(\Gamma)}.\eqno{(3.7)}$$
Let $\eta\in C_0^\infty(T^*\Gamma)$, $\eta=1$ on 
$\{\zeta\in T^*\Gamma:\|\zeta\|\le c_R^{-1}\}$,
$\eta=0$ on $\{\zeta\in T^*\Gamma:\|\zeta\|\ge c_R^{-1}+2\}$.
Since the parametrix ${\cal N}(\lambda)$ on supp$\,(1-\eta)$ is
an elliptic $\lambda-\Psi DO$ of class
$L^{1,0}_{0,0}(\Gamma)$, we have
$$\|{\rm Op}_\lambda(1-\eta)f\|_{H^1(\Gamma)}\le C
\left\|{\cal N}(\lambda)f\right\|_{L^2(\Gamma)}
+O(\lambda^{-\infty})\|f\|_{H^1(\Gamma)}$$ $$
\le C\|g\|_{L^2(\Gamma)}+C\|f\|_{L^2(\Gamma)}
+O(\lambda^{-\infty})\|f\|_{H^1(\Gamma)}.\eqno{(3.8)}$$
On the other hand,
$$\|{\rm Op}_\lambda(\eta)f\|_{H^1(\Gamma)}\le C
\|f\|_{L^2(\Gamma)}.\eqno{(3.9)}$$
Now (3.5) follows from combining (3.8) and (3.9) with (3.7).
\eproof

M. Khenissi,  D\'epartement de Math\'ematiques, Facult\'e des Sciences
de Monastir, Tunisie

e-mail: moez.khenissi@fsg.rnu.tn

G. Vodev, Universit\'e de Nantes,
 D\'epartement de Math\'ematiques, UMR 6629 du CNRS,
 2, rue de la Houssini\`ere, BP 92208, 44332 Nantes Cedex 03, France

e-mail: georgi.vodev@math.univ-nantes.fr

\end{document}